\documentclass[article,onesided]{amsart}\def\loadTIKZ{\usepackage{tikz}\usetikzlibrary{matrix,arrows,calc,cd,decorations.pathmorphing}}\overfullrule=5pt \newcommand\hfuzzReset{\hfuzz=3pt}\hfuzzReset \newcommand\toleranceReset{\tolerance=1400}\toleranceReset \newcommand\emergencystretchReset{\emergencystretch=1ex}\emergencystretchReset \hbadness=10000 \usepackage{ifpdf}\newcommand{\emailTressl}{marcus.tressl@manchester.ac.uk}\newcommand{\homepageTressl}{\url{http://personalpages.manchester.ac.uk/staff/Marcus.Tressl/}}\let\oldtocsection=\tocsection \let\oldtocsubsection=\tocsubsection \let\oldtocsubsubsection=\tocsubsubsection \renewcommand{\tocsection}[2]{\hspace{0em}\vspace*{0.1em}\oldtocsection{#1}{#2}}\renewcommand{\tocsubsection}[2]{\hspace{4ex}\oldtocsubsection{#1}{#2}}\renewcommand{\tocsubsubsection}[2]{\hspace{6ex}\oldtocsubsubsection{#1}{#2}}\ifpdf \usepackage[pdftex]{lscape}\else \usepackage{lscape}\fi \usepackage{ulem}\usepackage{fancybox}\usepackage{xifthen}\usepackage{forarray}\usepackage{xstring}\usepackage{stringstrings}\def\StackCreate#1#2#3{\expandafter\def\csname#1\endcsname{#3}\expandafter\def\csname#1Push\endcsname##1{\expandafter\edef\csname#1\endcsname{##1#2\csname#1\endcsname}}\expandafter\def\csname TopAux#1\endcsname ##1#2##2#3{##1}\expandafter\def\csname#1Top\endcsname{\expandafter\expandafter\expandafter\expandafter\expandafter\expandafter\csname TopAux#1\endcsname\csname#1\endcsname}\expandafter\def\csname PopAux#1\endcsname ##1#2##2#3##3#2{\expandafter\def\csname##3\endcsname{##2#3}}\expandafter\def\csname#1Pop\endcsname{\expandafter\expandafter\expandafter\expandafter\expandafter\expandafter\csname PopAux#1\endcsname\csname#1\endcsname#1#2}}\def\GetAfterColonAux#1:#2;{#2}\def\GetAfterColon#1{\IfBeginWith{#1}{:}{\GetAfterColonAux#1;}{#1}}\usepackage{aliascnt}\usepackage[shortlabels,inline]{enumitem}\setenumerate[1]{leftmargin=5.5ex}\setitemize[1]{leftmargin=5.5ex}\SetEnumitemKey{noindent}{leftmargin=0ex, itemindent=5ex, align=right, itemsep=1ex }\newcommand\NOPAGENUMBER[1]{}\usepackage{everypage}\newcommand\AddPrivateToMargin[1]{\AddEverypageHook{\tikz[overlay,remember picture]{\node at ($(current page.west)+(1.5,0)$) [rotate=90] {\textcolor{orange}{\vbox{\hrule width \the\textwidth height 0.5pt} \textcolor{black}{#1}\ \vbox{\hrule width 40em height 0.5pt}}}; }}}\newcommand\AddLongversionToMargin[1]{\AddEverypageHook{\tikz[overlay,remember picture]{\node at ($(current page.west)+(2,0)$) [rotate=90] {\textcolor{\LongColor}{\vbox{\hrule width \the\textwidth height 0.5pt} \textcolor{black}{#1}\ \vbox{\hrule width 40em height 0.5pt}}}; }}}\newcommand\AddOldversionToMargin[1]{\AddEverypageHook{\tikz[overlay,remember picture]{\node at ($(current page.west)+(2.5,0)$) [rotate=90] {\textcolor{\OldColor}{\vbox{\hrule width \the\textwidth height 0.5pt} \textcolor{black}{#1}\ \vbox{\hrule width 40em height 0.5pt}}}; }}}\newcommand\AddLineToMargin[3]{\AddEverypageHook{\tikz[overlay,remember picture]{\node at ($(current page.west)+(#2,0)$) [rotate=90] {\textcolor{#1}{\vbox{\hrule width \the\textwidth height 0.5pt} \textcolor{black}{#3}\ \vbox{\hrule width 40em height 0.5pt}}}; }}}\IfFileExists{mathabx.sty}{}{}\usepackage{amsfonts}\usepackage{amssymb}\usepackage{stmaryrd}\usepackage{amsmath}\usepackage{amsthm}\usepackage{dsfont}\IfFileExists{mbboard.sty}{\usepackage{mbboard}}{}\usepackage{mathrsfs}\usepackage{twcal}\usepackage{accents}\usepackage{bm}\usepackage[T1]{fontenc}\usepackage[latin1]{inputenc}\ifpdf \usepackage[pdftex,usenames,x11names]{xcolor}\else \usepackage[dvips,usenames,x11names]{xcolor}\fi \usepackage[pdftex]{graphicx}\usepackage[all]{xy}\ifdefined\loadTIKZ \loadTIKZ \def\TIKZlabel#1{}\else\fi \StackCreate{ColoR}{;}{?}\ColoRPush{black}\newcommand {\notion}[2][]{\def\temp{#1}\ifmmode #2 \ifx \temp\empty \index{$#2$}\else \index{$#1$}\fi \else {\bf #2}\ifx \temp\empty \index{#2}\else \index{#1}\fi \fi }\ifpdf \usepackage[pdftex,linktocpage,pagebackref,breaklinks]{hyperref}\hypersetup{colorlinks=true,allcolors=Green, linkcolor=DarkGreen, citecolor=violet, urlcolor=blue, runcolor=red, filecolor=cyan }\else \usepackage[hypertex,linktocpage,pagebackref]{hyperref}\fi \def\UndefinedRef#1{\LARGE\bfseries\color{red} ??#1??}\makeatletter \def\@setref#1#2#3{\ifx#1\relax \protect\G@refundefinedtrue \nfss@text{\reset@font\UndefinedRef{#3}}\@latex@warning{Reference `#3' on page \thepage \space undefined }\else \expandafter\Hy@setref@link#1\@empty\@empty\@nil{#2}\fi }\makeatother \usepackage{xr-hyper}\newcommand{\refX}[2]{\IfBeginWith{#1}{:}{\ref{\GetAfterColonAux#1;-#2}}{\cite[\ref{#1-#2}]{#1}}}\newcommand\pr{\begin{proof}}\def\ende{\end{proof}}\newtheoremstyle{LayoutVoid}{1ex}{0ex}{\normalfont}{}{\bf}{.}{1ex}{}\newcommand\stressstatement[1]{#1}\theoremstyle{plain}\swapnumbers \newcommand\maketheorem[1]{\newtheorem{#1}[theorem]{\stressstatement{#1}} \newtheorem{#1Definition}[theorem]{\stressstatement{#1 and Definition}}  }\FunctionForEach{,}{\maketheorem}{Conclusion,Conjecture,Corollary,Fact,Facts,Lemma,Observation,Observations,Proposition,Reminder,Scholium,Summary,Theorem}\theoremstyle{definition}\theoremstyle{remark}\FunctionForEach{,}{\maketheorem}{Convention,Counterexample,Discussion,Example,Examples, Exercise,Exercises,Explanation,Notation,Project,Projects,Question,Questions,Remark,Remarks,Strategy,Warning}\theoremstyle{LayoutVoid}\numberwithin{equation}{section}\newcommand{\labelon}[1]{\marginpar{#1}}\newcommand{\labelx}[1]{{\def\temp{#1}\ifx\temp\empty\else \label{#1}\labelon{#1}\fi}}\def\GetAfterColon#1:#2;;{#2}\def\GetAfterPlus#1+#2;;{#2}\newenvironment{FACT}[2]{\IfBeginWith{#1}{:}{\def\tempFactName{void}\def\tempFreeTitle{\GetAfterColon#1;;\ }}{\IfBeginWith{#1}{+}{\def\tempFactName{voidTheorem}\def\tempFreeTitle{\GetAfterPlus#1;;\ }}{\def\tempFactName{#1}\def\tempFreeTitle{}}}\def\tempfacT{\end{\tempFactName}}\begin{\tempFactName}\labelx{#2}\textup{\textbf{\tempFreeTitle}}\capitalize[q]{#1}\caselower[q]{#1}\global\edef\factname{\thestring}}{\tempfacT}\catcode`\=13 \def{+}\newcommand\assigncharacter[1]{\expandafter\newcommand\csname #1\endcsname{\mathds{#1}}}\FunctionForEach{,}{\assigncharacter}{A,B,C,D,E,F,G,I,J,K,M,N,Q,R,T,U,V,W,X,Y,Z}\renewcommand\assigncharacter[1]{\expandafter\newcommand\csname C#1\endcsname{\mathcal{#1}}}\FunctionForEach{,}{\assigncharacter}{A,B,C,D,E,F,G,H,I,J,K,L,M,N,O,P,Q,R,S,T,U,V,W,X,Y,Z}\renewcommand\assigncharacter[1]{\expandafter\newcommand\csname D#1\endcsname{\mathfrak{#1}}}\FunctionForEach{,}{\assigncharacter}{a,b,c,d,e,f,g,h,i,j,k,l,m,n,o,p,q,r,s,t,u,v,w,x,y,z,A,B,C,D,E,F,G,I,K,L,M,N,O,P,Q,R,S,T,U,V,W,X,Y,Z} \renewcommand\assigncharacter[1]{\expandafter\newcommand\csname S#1\endcsname{\mathscr{#1}}}\FunctionForEach{,}{\assigncharacter}{A,B,C,D,E,F,G,H,I,J,K,L,M,N,O,P,Q,R,T,U,V,W,X,Y,Z}\def\NewFont#1#2#3#4#5{\expandafter\font\csname #1display\endcsname =#1 at #2 \expandafter\font\csname #1normal\endcsname =#1 at #3 \expandafter\font\csname #1script\endcsname =#1 at #4 \expandafter\font\csname #1scriptscript\endcsname =#1 at #5 }\def\NewFontLetter#1#2{{\mathchoice {{\expandafter\hbox{\csname #1display\endcsname\char"#2}}}{{\expandafter\hbox{\csname #1normal\endcsname\char"#2}}}{{\expandafter\hbox{\csname #1script\endcsname\char"#2}}}{{\expandafter\hbox{\csname #1scriptscript\endcsname\char"#2}}}}}\NewFont{pxsyc}{9.00pt}{8.00pt}{7.00pt}{6.00pt}\NewFont{pxsya}{9.00pt}{8.00pt}{7.00pt}{6.00pt}\NewFont{p1xr}{10.00pt}{9.00pt}{8.00pt}{7.00pt}\NewFont{MnSymbolC10}{10.00pt}{9.00pt}{8.00pt}{7.00pt}\NewFont{MnSymbolD10}{12.00pt}{11.00pt}{10.00pt}{9.00pt}\NewFont{MnSymbolF10}{12.00pt}{11.00pt}{10.00pt}{9.00pt}\NewFont{manfnt}{12.00pt}{11.00pt}{10.00pt}{9.00pt}\NewFont{favmr7y}{12.00pt}{11.00pt}{10.00pt}{9.00pt}\catcode95=12 \catcode95=8 \newcommand\bdl{{\ifmmode \mathrm{bdlat}\else {bounded distributive lattice}\fi}} \newcommand{\st}{{\ \vert \ }}\let\temp\phi \let\phi\varphi \let \varphi\temp \let\temp\theta \let\theta\vartheta \let \vartheta\temp \let\eps\varepsilon  \let\0\emptyset \newcommand{\onto}{\twoheadrightarrow}\newcommand{\lra}{\longrightarrow}\newcommand{\mal}{\cdot}\newcommand{\monthname}[1]{\ifcase#1 \or January \or February \or March \or April \or May \or June \or July \or August \or September \or October \or November \or December \fi}\newcommand\LongColor{teal}\newcommand\OldColor{gray}\newcommand\COL{\ifmmode\colon\else :\ \fi}\newcommand\kat[1]{{\tt #1}}\newcommand{\Claim}[1]{\underline{Claim #1.}\ }\newcommand\map[5]{\begin{eqnarray*}#1#2&\lra &#3 \\ #4&\longmapsto &#5 \end{eqnarray*}} \renewcommand{\mod}{{\operator{\,mod\,}}}\newcommand\operator[1]{\mathop{\operatorname{#1}}\nolimits} \newcommand{\id}{\operator{id}}\newcommand{\Der}{\operator{Der}}\newcommand{\Hom}{\operator{Hom}}\newcommand\DIRroot{../../../../../../../../../../../}\definecolor{Green}{rgb}{0.00,0.50,0.00}\definecolor{DarkGreen}{rgb}{0.00,0.40,0.00}\definecolor{grey}{rgb}{0.40,0.40,0.40} \renewcommand\textcolor[2]{\ColoRPush{#1}\color{\ColoRTop}#2\ColoRPop\color{\ColoRTop}}\IfFileExists{C:/wb/System64/WinBatch.exe}{}{}\let\temp\theta \let\theta\vartheta \let \vartheta\temp \def\Ind{\setbox0=\hbox{$x$}\kern\wd0\hbox to 0pt{\hss$\mid$\hss} \lower.9\ht0\hbox to 0pt{\hss$\smile$\hss}\kern\wd0}\def\Notind{\setbox0=\hbox{$x$}\kern\wd0\hbox to 0pt{\mathchardef \nn=12854\hss$\nn$\kern1.4\wd0\hss}\hbox to 0pt{\hss$\mid$\hss}\lower.9\ht0 \hbox to 0pt{\hss$\smile$\hss}\kern\wd0}\renewcommand{\labelon}[1]{}
\begin{document} \title{Differential Weil descent} \author{Omar Le\'on S\'anchez} \address{Omar Le\'on S\'anchez, The University of Manchester\\ Department of Mathematics\\ Oxford Road \\ Manchester, M13 9PL, UK} \email{omar.sanchez@manchester.ac.uk} \author{Marcus Tressl} \address{Marcus Tressl, The University of Manchester\\ Department of Mathematics\\ Oxford Road \\ Manchester, M13 9PL, UK \newline Homepage: \homepageTressl} \email{\emailTressl} \date{\today} \subjclass[2010]{12H05, 14A99} \keywords{differential algebras, Weil descent} \begin{abstract} In this short note a differential version of the classical Weil descent is established in all characteristics. This yields a ready-to-deploy tool of differential restriction of scalars for differential varieties over finite differential field extensions. \end{abstract} \maketitle \section{Introduction} \noindent In 1959, Andr\'e Weil introduced his method of restricting scalars for a finite separable field extension $K\subseteq L$, cf. \cite[\S1.3]{Weil1982}. It says that scalar extension, seen as a functor from $K$-algebras to $L$-algebras, has a left adjoint, which sends an $L$-algebra $D$ to a $K$-algebra $W(D)$, the Weil descent (aka Weil restriction) of $D$ from $L$ to $K$. The construction has been vastly generalised by Grothendieck \cite{Grothe1995a}, and used in numerous occasions in number theory \cite{SmartWebsite} and algebraic geometry \cite{Milne1972}. \par We establish a similar descent for differential algebras with respect to a given extension of differential rings $A\subseteq B$, where $B$ is finitely generated and free as an $A$-module. Here a differential ring $A$ is a commutative unital ring equipped with a distinguished set of derivations $A\lra A$. If $D$ is a differential $B$-algebra with commuting derivations, its descent $W^\mathrm{diff}(D)$ is a differential $A$-algebra in commuting derivations, see Theorem~\ref{DiffWeilDescent}. This is deduced from our main result, which concerns rings and algebras with a single derivation: \par \medskip\noindent \textbf{Main Theorem} (see Theorems \ref{partialDW} and \ref{liebracket}) \par \noindent Let $d:A\lra A$ be a derivation of a ring $A$ and let $(B,\delta)$ be a differential $(A,d)$-algebra. Assume that $B$ is finitely generated and free as an $A$-module. \begin{enumerate}[(i),itemsep=1ex] \item Let $(D,\partial)$ be a differential $(B,\delta)$-algebra. Then there is a unique derivation $\partial^W$ on the classical Weil descent $W(D)$ such that $(W(D),\partial^W)$ is a differential $(A,d)$-algebra and the unit of the adjunction at $D$ (given by the classical Weil descent), namely the map $W_D:D\lra W(D)\otimes B$, is a differential $(B,\delta)$-algebra homomorphism $ (D,\partial)\lra (W(D)\otimes B,\partial^W\otimes \delta). $ \item If $B$ is a subring of $D$ and the inclusion is the structure morphism of $D$ as a $B$-algebra, then the assignment $\partial\mapsto \partial^W$ is a Lie-ring and an $A$-module homomorphism. \end{enumerate} \par \noindent A special case of this theorem appears in \cite[\S5]{LeSMos2016}, when $A$ and $B$ are differential fields (of characteristic 0), but only under the assumption that $B$ has an $A$-basis consisting of constants, see \ref{remarkOnLeSMos}. However, a basis of constants does not always exist, as we point out in Example \ref{example1}. \par \smallskip In the proof of the main theorem, we give explicit formulas of the involved differential rings and morphisms, rather than only showing that they exist. \par \smallskip In a forthcoming paper by the authors several applications of the differential Weil descent are exposed. The main one addresses a method to produce differential fields, in finitely many commuting derivations and of characteristic 0, which possess a minimal differential closure (or in Kolchin's terminology, \textit{constraint closure}). First examples of such differential fields were given by Singer in \cite{Singer1978b}, where he showed that, for every closed ordered differential field $K$ in one derivation, the algebraic closure $K[i]$ is differentially closed. In the coming paper, we introduce the notion of a differentially large field (in analogy to the notion of a large field in classical field theory) and use the differential Weil descent to generalize Singer's result; namely, we will show that algebraic extensions (equipped with the unique induced derivations) of differentially large fields are again differentially large. This is in analogy to the algebraic case where the classical Weil descent is used to show that algebraic extensions of large fields are again large, see \cite[Proposition 1.2]{Pop1996}. \par \smallskip We expect many more applications of the differential Weil descent. For instance, we expect that our results will be a valuable tool in differential Galois cohomology and the parameterised Picard-Vessiot theory for linear differential equations. This will potentially be in the form of finiteness results for cohomology groups of linear differential algebraic groups. \section{Classical Weil Descent for Algebras}\labelx{classicalweil} \labelx{ClassicalWeil} \noindent In this section we review the classical construction of Weil descent of scalars for algebras, see for example \cite[\S7.6]{BoLuRa1990}, \cite[\S2]{MooSca2010} and \cite{Grothe1995a}. For our purposes we need certain explicit formulas, so we give details. \par \smallskip\noindent \textbf{Convention.} Throughout, we assume our rings and algebras to be commutative and unital; ring and algebra homomorphisms are meant to be unital as well. \par \smallskip Let $A$ be a ring and let $B$ be an $A$-algebra. For each $A$-algebra $C$, the scalar extension by $B$ is the $B$-algebra $C\otimes_AB$ with structure map $b\mapsto 1\otimes b$ \footnote{As a general reference for tensor products, specifically in the category of algebras we refer to \cite[Appendix A]{Matsum1989}}. This assignment has a natural extension to a covariant functor $F:A\text{-}\kat{Alg}\lra B\text{-}\kat{Alg}$. The functor $F$ has a right adjoint $B\text{-}\kat{Alg}\lra A\text{-}\kat{Alg}$ given by restricting scalars. If $B$ is finitely generated and free as an $A$-module, then $F$ also has a left adjoint $W$, called \notion{Weil descent}, or \notion{Weil restriction}. We start with a reminder on left adjoints in general, ready made for use later on. \begin{FACT}{:Fact.}{LeftAdjoint} \textup{\cite[Thm 2, p.83, Cor. 1,2, p.84]{MacLan1998}} \noindent Let $F:\CC\lra \CD$ be a covariant functor between categories $\CC$ and $\CD$. \begin{enumerate}[(i)] \item The following are equivalent. \begin{enumerate} \item\labelx{LeftAdjointLA}$F$ has a left adjoint $W$, i.e., $W:\CD\lra \CC$ is a covariant functor such that for all $D\in \CD$ the functor $\Hom_\CD(D,F(\ \underline {\ }\ )):\CC\lra \kat{Sets}$ is represented by $W(D)$, meaning that the functors $\Hom_{{\CC}}({W}(D),\ \underline {\ }\ )$ and $\Hom_\CD(D,F(\ \underline {\ }\ ))$ are isomorphic \footnote{Recall that two functors are isomorphic if there is an invertible natural transformation between them.}. \item\labelx{LeftAdjointDagger} For each $D\in \CD$ there are $W(D)\in \CC$ and a $\CD$-morphism $W_D:D\lra F(W(D))$ such that the following condition holds: \par \smallskip \begin{enumerate} \item[$(\dagger)$:] For every $C\in\CC$ and each morphism $f:D\lra F(C)$, there is a unique $\CC$-morphism $g:{W}(D)\lra C$ such that the following diagram commutes \begin{center} \begin{tikzcd} F(W(D)) \ar[r, dashed, "F(g)"] & F(C)\\ D \ar[u, "W_D"] \ar[ru,"f"'] \end{tikzcd} \end{center} \noindent In other words, $W_D$ gives rise to a bijection \[ \tau(D,C):\Hom_{{\CC}}({W}(D),C)\lra \Hom_{{\CD}}(D,{F}(C)),\ g\longmapsto {F}(g)\circ {W}_D. \] \end{enumerate} \end{enumerate} \end{enumerate} \begin{enumerate}[(i),resume] \item\labelx{LeftAdjointLAimpliesDagger} If \ref{LeftAdjointLA} holds, then for every such functor $W$, all $D\in\CD$ and each isomorphism $\tau(D,\ \underline{\ }\ ) :\Hom_{{\CC}}({W}(D),\ \underline {\ }\ )\lra \Hom_{\CD}(D,F(\ \underline {\ }\ ))$ as in \ref{LeftAdjointLA}, the choice $W(D)$ and $W_D=\tau(D,W(D))(\id_{W(D)})$ satisfy property $(\dagger)$ of \ref{LeftAdjointDagger}. \par The assignment $D\mapsto W_D$ is a natural transformation $\id_{\CD}\lra F\circ W$ and is called the \notion{unit of the adjunction}; $W_D$ is called the \notion{component} at $D$ of that unit. \par Similarly, for each $C\in \CC$ the morphism $F_C:W(F(C))\lra C$ that is sent to $\id_{F(C)}$ via $\tau(F(C),C)$ gives rise to a natural transformation $W\circ F\lra \id_{\CC}$, called the \notion{counit of the adjunction}; $F_C$ is called the \notion{component} at $C$ of that counit. \item\labelx{LeftAdjointDaggerImpliesLA} If \ref{LeftAdjointDagger} holds, then for any choice of objects $W(D)$ and morphisms $W_D$ as in \ref{LeftAdjointDagger}, for $D\in\CD$, the assignment $D\mapsto W(D)$ can be extended to a functor $W:\CD\lra\CC$ satisfying \ref{LeftAdjointLA} as follows: Take a morphism $f_0:D\lra D'$ and set $f:=W_{D'}\circ f_0$ and $C:=W(D')$. Then define $W(f_0):W(D)\lra W(D')$ as the unique $\CC$-morphism $W(D)\lra W(D')$ such that the diagram \begin{center} \begin{tikzcd}[row sep=7ex,column sep=10ex] F(W(D)) \ar[r, dashed, "F(W(f_0))"] & F(C)=F(W(D'))\\ D \ar[u, "W_D"]\ar[r,"f_0"'] \ar[ru,"f"'] & D'\ar[u, "W_{D'}"] \end{tikzcd} \end{center} commutes, according to $(\dagger)$. \item\labelx{LeftAdjointUnique} Any two functors that are left adjoint to $F$ are isomorphic. \item\labelx{LeftAdjointExact} If $W$ is left adjoint to $F$, then $W$ preserves all co-limits, cf. \cite[p. 119, last paragraph]{MacLan1998}. For example $W$ preserves direct limits and fiber sums (aka pushouts). \par \end{enumerate} \end{FACT} \begin{FACT}{:Notation and setup.}{NotationSetup} We return to our setup of a ring $A$ and an $A$-algebra $B$. Let \[ F:A\text{-}\kat{Alg}\lra B\text{-}\kat{Alg} \] be the functor defined by $F(C)=C\otimes B$ and for $\phi:C\lra C'$, $F(\phi)=\phi\otimes \id_B$. Here and below, tensor products are taken over $A$, unless stated otherwise. \par \smallskip\noindent We will from now on assume that $B$ is free and finitely generated as an $A$-module of dimension $\ell$ over $A$. We also fix generators $b_1,\ldots,b_\ell $ of the $A$-module $B$. For $i\in \{ 1,...,\ell \} $ let \[ \lambda_i:B\lra A,\ \lambda_i (\sum_{j=1}^\ell a_jb_j)=a_i \] be the $A$-module homomorphism dual to $A\lra B,\ a\mapsto a\mal b_i$. If $C$ is an $A$-algebra we write $\lambda_i^C=\id_C\otimes \lambda_i:C\otimes B\lra C\otimes A=C$ for the base change of $\lambda_i$ to $C$. Since $b=\sum_i \lambda_i(b)b_i$ for $b\in B$ we obtain \begin{align*} c\otimes b=c\otimes \sum_i \lambda_i(b)b_i =\sum_i (\lambda_i(b)c)\otimes b_i, \end{align*} for $c\in C$. Hence $\lambda_i^C(c\otimes b)=\lambda_i(b)\mal c$ is the coefficient of $c\otimes b$ at $1\otimes b_i$ when it is written in the basis $1\otimes b_1,\ldots,1\otimes b_\ell $ of the free $C$-module $C\otimes B$. This extends to all $f\in C\otimes B$, thus \[ f=\sum_{i=1}^\ell \bigl(\lambda_i^C(f)\otimes b_i\bigr). \leqno{(*)} \] \end{FACT} \begin{FACT}{Definition}{DefnWBT} Let $T$ be a set of indeterminates for $A$ and $B$. We define an $A$-algebra $W(B[T])=A[T]^{\otimes \ell}\ (=\underbrace{A[T]\otimes\ldots\otimes A[T]}_{\ell\text{-times}})$. For $i\in \{1,\ldots,\ell \}$ and $t\in T$ we write \[ t(i):=1\otimes ...\otimes 1\otimes \underbrace{t}_{i\text{-th position}}\otimes 1\otimes ...\otimes 1\in A[T]^{\otimes \ell}. \] Let $W_{B[T]}$ be the unique $B$-algebra homomorphism \[ W_{B[T]}:B[T]\lra F(W(B[T]))=A[T]^{\otimes \ell}\otimes B\text{ with } \] \[ W_{B[T]}(t)=\sum _{i=1}^\ell (t(i)\otimes b_i)\quad (t\in T). \] Further, let $F_{A[T]}$ be the unique $A$-algebra homomorphism \[ F_{A[T]}:W(F(A[T]))=A[T]^{\otimes \ell}\lra A[T] \] with the property $F_{A[T]}(t(i))=\lambda_i(1)\mal t$ for $t\in T,\ i\in\{1,\ldots,\ell \}$. \end{FACT} \begin{FACT}{:Explicit description of the Weil descent of polynomial algebras.}{WeilBT} The $A$-algebra $W(B[T])$ and the morphism $W_{B[T]}$ described above satisfy condition $(\dagger)$ of \ref{LeftAdjoint}\ref{LeftAdjointDagger}. Hence by \ref{LeftAdjoint}\ref{LeftAdjointDaggerImpliesLA} we may choose $W(B[T])$ as the Weil descent of $B[T]$, and $W_{B[T]}$ as the unit of the adjunction at $B[T]$; these choices are then independent of the basis $b_1,\ldots,b_\ell$ up to a natural $A$-algebra isomorphism (see \ref{LeftAdjoint}\ref{LeftAdjointUnique}). \par Explicitly, for every $C\in A\text{-}\kat{Alg}$, the map \map{\tau=\tau(B[T],C):\ }{\Hom_{A\text{-}\kat{Alg}}(A[T]^{\otimes \ell},C)} {\Hom_{B\text{-}\kat{Alg}}(B[T],C\otimes B)}{\phi} {F(\phi)\circ W_{B[T]}=(\phi \otimes \id_B)\circ W_{B[T]}} is bijective, where $\phi \otimes \id_B=F(\phi):F(W(B[T]))=A[T]^{\otimes d}\otimes B\lra C\otimes B$ is the base change of $\phi $. For $t\in T$ we have \[ \tau(\phi)(t)=\sum _{i=1}^\ell (\phi(t(i))\otimes b_i). \] \par \smallskip\noindent The compositional inverse of $\tau=\tau(B[T],C)$ is defined as follows. Let $\psi :B[T]\lra C\otimes B$ be a $B$-algebra homomorphism. We define an $A$-algebra homomorphism $\phi :A[T]^{\otimes d}\lra C$ by \[ \phi (t(i)):=\lambda_i^C(\psi(t))\ (t\in T,\ i=1,\ldots,\ell ). \] Since $A[T]\otimes B\cong _BB[T]$, $\psi $ is uniquely determined by $\{\psi(t)\st t\in T\}$ and we see that $\phi $ is the unique preimage of $\psi $ under $\tau$. \par \noindent Further, one checks easily that $F_{A[T]}$ is the component of the counit of the adjunction at $A[T]$. \end{FACT} \begin{FACT}{:Explicit description of the Weil descent of $B$-algebras.}{WeilD} Now let $D$ be a $B$-algebra. Take a surjective $B$-algebra homomorphism $\pi_D:B[T]\onto D$ for some set $T$ of indeterminates. Let $I_D$ be the ideal generated in $W(B[T])=A[T]^{\otimes \ell}$ generated by all the $\lambda_i^{W(B[T])}(W_{B[T]}(f))$, where $i\in\{1,\dots,\ell\}$ and $f\in\ker (\pi_D)$. We define \[ W(D):=W(B[T])/I_D \] and write \[ W(\pi_D):W(B[T])\lra W(D) \] for the residue map. Then the bijection $\tau(B{[}T{]}{,}C)$ from \ref{WeilBT} induces a bijection \[ \tau(D{,}C):\Hom_{A\text{-}\kat{Alg}}(W(D),C)\lra \Hom_{B\text{-}\kat{Alg}}(D,F(C)) \] such that the diagram \begin{center} \begin{tikzcd}[row sep=10ex,column sep=10ex] \Hom_{A\text{-}\kat{Alg}}(W(D),C) \ar[d,"\underline{\ }\circ W(\pi_D)"', hook] \ar[r,"\tau(D{,}C)", dashed] & \Hom_{B\text{-}\kat{Alg}}(D,F(C))\ar[d,"\ \underline{\ }\circ \pi_D", hook] \\ \Hom_{A\text{-}\kat{Alg}}(W(B[T]),C) \ar[r,"\tau(B{[}T{]}{,}C)"', "\simeq"] & \Hom_{B\text{-}\kat{Alg}}(B[T],F(C)) \end{tikzcd} \end{center} commutes. \par \noindent The commutativity of the diagram above says that for $\phi\in \Hom_{A\text{-}\kat{Alg}}(W(D),C)$ we have \begin{align*} (+)\quad \tau(D,C)(\phi)\circ \pi_D&=\tau(B{[}T{]}{,}C)(\phi\circ W(\pi_D))\cr &=((\phi\circ W(\pi_D))\otimes\id_B )\circ W_{B[T]}. \end{align*} \par \noindent Finally, we display the map $W_D:=\tau(D,W(D))(\id_{W(D)}):D\lra F(W(D))$ explicitly and show that -- together with $W(D)$ -- it satisfies the mapping property of $(\dagger)$ in \ref{LeftAdjoint}\ref{LeftAdjointDagger}. Take $t\in T$. Then by (+) with $C=W(D),\ \phi=\id_{W(D)}$ we see that \[ W_D(\pi_D(t))=\sum _{i=1}^\ell W(\pi_D)(t(i))\otimes b_i=\sum _{i=1}^\ell (t(i)\mod I_D)\otimes b_i.\leqno{(\sharp)} \] \par \noindent Pick an $A$-algebra $C$. Since $\tau(D,C)$ is bijective, the mapping property of $(\dagger)$ in \ref{LeftAdjoint}\ref{LeftAdjointDagger} follows after checking $\tau(D,C)(\phi)=F(\phi)\circ W_D$ for all $\phi\in \Hom_{A\text{-}\kat{Alg}}(W(D),C)$. Using $(\sharp)$ this is a straightforward computation. \par Using \ref{LeftAdjoint}\ref{LeftAdjointDaggerImpliesLA},\ref{LeftAdjointUnique} we have justified our choice of $W(D)$ and $W_D$ for the Weil descent. Finally, \ref{LeftAdjoint}\ref{LeftAdjointDaggerImpliesLA} gives the definition of $W$ on morphisms. \end{FACT} \section{Differential Weil Descent}\label{differentialweil} \noindent In this section we present a construction of a Weil descent functor in the category of differential algebras in arbitrary characteristic. We first recall some basic facts about differential algebras and their tensor products. We continue to assume that our rings and algebras are unital and commutative. \begin{FACT}{:Generalities about differential algebra.}{GenDiff} The following are well known generalities on differential algebras whose proofs are straightforward. For a ring $A$ we let $\Der(A)$ denote the family of derivations on $A$. \begin{enumerate}[(i)] \item\labelx{GenDiffExtendOrdinary} Let $A$ be a ring and let $T$ be a not necessarily finite set of indeterminates over $A$. For each $t\in T$ let $f_t\in A[T]$. Let $d\in \Der(A)$. Then there is a unique derivation $\delta $ of $A[T]$ extending $d$ with $\delta(t)=f_t$ for all $t\in T$. \par \end{enumerate} \par \noindent For $ d,\delta\in\Der(A)$ we write $[ d,\delta]:A\lra A$ for the Lie-bracket of $ d $ and $\delta$, defined by $[ d,\delta](a)= d\delta(a)-\delta d(a)$. Notice that $[ d,\delta]$ is again a derivation of $A$. \begin{enumerate}[(i),resume] \item \labelx{GenDiffLieFromGenerators} Let $A$ be a ring and let $S\subseteq A$ be a set of generators of the ring $A$. \begin{enumerate}[(a)] \item Let $ d$ and $(\delta_i)_{i\in I}$ be derivations on $A$ and suppose there are $a_{i}\in A$, all but finitely many zero, with $d(s)=\sum_{i\in I}a_{i}\delta_i(s)$ for all $s\in S$. Then $d=\sum_{i\in I}a_{i}\delta_i$. \item Let $\phi:A\lra B$ be a ring homomorphism and let $ d:A\lra A,\ \delta:B\lra B$ be derivations. If $\phi( d s)=\delta(\phi(s))$ for all $s\in S$, then $\phi$ is a differential homomorphism $(A, d)\lra (B,\delta)$. \par \end{enumerate} \item \labelx{GenDiffTensor} Let $d\in\Der(A)$ and let $(B,\delta), (C,\partial)$ be differential $(A,d)$-algebras. Then there is a unique derivation $\delta\otimes \partial$ on $B\otimes_AC$ such that the natural maps $B\lra B\otimes_AC, C\lra B\otimes_AC$ are differential maps, cf. \cite[Chapter 2 (1.1), p. 21]{Buium1994}. \item\labelx{GenDiffLieAndTensor} Now let $d_1,d_2\in\Der(A)$, $\delta_1,\delta_2\in\Der(B)$ and $\partial_1,\partial_2\in\Der(C)$ such that $(B,\delta_i),(C,\partial_i)$ are differential $(A,d_i)$-algebras. Then, for $a_1,a_2\in A$, straightforward checking shows that \begin{enumerate} \item $(a_1\delta_1+a_2\delta_2)\otimes (a_1\partial_1+a_2\partial_2)=a_1(\delta_1\otimes \partial_1)+ a_2(\delta_2\otimes \partial_2)$. \item $[\delta_1,\delta_2]\otimes [\partial_1,\partial_2]=[\delta_1\otimes \partial_1,\delta_2\otimes \partial_2].$ \end{enumerate} \end{enumerate} \end{FACT} \bigskip \par \noindent As in Section \ref{classicalweil} we work with a ring $A$ and an $A$-algebra $B$ that is free and finitely generated by $b_1,\ldots,b_\ell $ as an $A$-module. We fix a derivation $d$ on $A$ and a derivation $\delta$ on $B$ such that $(B,\delta)$ is a differential $(A,d)$-algebra (meaning that the structure map $A\lra B$ is differential). \par By \ref{GenDiff}\ref{GenDiffTensor}, for any differential $(A,d)$-algebra $(C,\partial_C)$, there is a unique derivation $\partial_C\otimes \delta$ on $F(C)=C\otimes B$ such that the natural map $C\lra F(C)$ is a differential $(B,\delta)$-algebra morphism. \begin{FACT}{Theorem}{partialDW} Let $(D,\partial_D)$ be a differential $(B,\delta)$-algebra. Then there is a unique derivation $\partial_D^W$ on $W(D)$ such that $(W(D),\partial_D^W)$ is a differential $(A,d)$-algebra and \[ W_D:(D,\partial_D)\lra (F(W(D)),\partial_D^W\otimes \delta) \] is a differential $(B,\delta)$-algebra homomorphism, i.e., $W_D\circ \partial_D=(\partial_D^W\otimes \delta)\circ W_D$. \par Furthermore, $\partial_D^W$ only depends on $\partial_D$ and not on $\delta$. \end{FACT} \begin{proof} Take any set $T$ of differential indeterminates and a surjective $(B,\delta)$-algebra homomorphism $\pi_D:(B\{T\},\partial)\onto (D,\partial_D)$. Here, the differential polynomial ring $B\{T\}$ is considered just as polynomial ring over $B$ in the algebraic indeterminates $t_\theta$, where $t\in T$ and $\theta\in\Theta:=\{\partial^i:i\geq 0\}$. Further, $\partial=\partial_{B\{T\}}:B\{T\}\lra B\{T\}$ is the natural derivation, thus $\partial t_\theta =t_{\partial\theta}$. \par We choose $W_{B\{T\}}:B\{T\}\lra F(W(B\{T\}))$ according to \ref{DefnWBT} for the set of indeterminates $\{t_\theta\st t\in T,\ \theta\in\Theta\}$ and $W_D:D\lra F(W(D))$ according to \ref{WeilD}. Also recall $(\sharp)$ in \ref{WeilD}, which says that $W_D(\pi_D(t_\theta))=\sum _{i=1}^\ell W(\pi_D)(t_\theta(i))\otimes b_i$ \par \smallskip \par \noindent \Claim 1 If $\eps:W(D)\lra W(D)$ is a derivation such that $(W(D),\eps)$ is a differential $A$-algebra, then for all $t\in T$ and any $\theta\in \Theta$ we have \[ ((\eps\otimes \delta)\circ W_D)(\pi_D(t_\theta))= \sum_{i=1}^\ell \biggl(\eps (W(\pi_D)(t_\theta(i)))+\sum_{j=1}^\ell \lambda_i(\delta b_j)\mal W(\pi_D)(t_\theta(j)) \biggr)\otimes b_i. \] See \ref{GenDiff}\ref{GenDiffTensor} for the definition of $\eps\otimes\delta$. \par \noindent \textit{Proof.} This is a straightforward calculation using $\delta b_i=\sum_{j=1}^\ell \lambda_j(\delta b_i)b_j$. \hfill$\diamond$ \par \smallskip\noindent \Claim 2 If $\eps:W(D)\lra W(D)$ is a derivation such that $(W(D),\eps)$ is a differential $A$-algebra, then $W_D\circ \partial=(\eps\otimes \delta)\circ W_D$ if and only if for all $t_\theta(i)$ we have \[ \eps(W(\pi_D)(t_\theta(i)))= W(\pi_D)(t_{\partial\theta}(i)-\sum_{j=1}^\ell \lambda_i(\delta(b_j))\mal W(\pi_D)(t_\theta(j)).\leqno{(*)} \] \noindent \textit{Proof.} By \ref{GenDiff}\ref{GenDiffLieFromGenerators}(b), $W_D\circ \partial=(\eps\otimes \delta)\circ W_D$ if and only if $((\eps\otimes \delta)\circ W_D)(\pi_D(t_\theta))=(W_D\circ \partial)(\pi_D(t_\theta))$ for all $t_\theta$. By Claim 1 this is equivalent to \begin{align*} \sum_{i=1}^\ell &\biggl(\eps (W(\pi_D)(t_\theta(i)))+\sum_{j=1}^\ell \lambda_i(\delta b_j)\mal W(\pi_D)(t_\theta(j)) \biggr)\otimes b_i\cr &=W_D(\partial(\pi_D(t_\theta)))\cr &=W_D(\pi_D(\partial t_\theta)),\text{ since }\pi_D\text{ is a differential map}\cr &=W_D(\pi_D(t_{\partial\theta}))\cr &=\sum _{i=1}^\ell W(\pi_D)(t_{\partial\theta}(i))\otimes b_i\text{, by }(\sharp)\text{ in \ref{WeilD}}. \end{align*} Since $1\otimes b_1,\ldots,1\otimes b_\ell$ is a basis of $F(W(D))$ over $W(D)$, the identity is equivalent to $(*)$ being true for all $i\in\{1,\ldots,\ell\}$. \hfill$\diamond$ \par \medskip\noindent Claim 2 implies the uniqueness statement of the \factname, because the set of all the $W(\pi_D)(t_\theta(i))$ generates $W(D)$. For existence, we first deal with $B\{T\}$ instead of $D$. In that case, Claim 2 says that we only need to find a derivation $\partial_{B\{T\}}^W$ on $W(B\{T\})$ such that $(W(B\{T\}),\partial_{B\{T\}}^W)$ is a differential $(A,d)$-algebra with the property \[ \partial_{B\{T\}}^W(t_\theta(i))=t_{\partial \theta}(i)-\sum_{j=1}^\ell \lambda_i(\delta(b_j))\mal t_\theta(j). \] By \ref{GenDiff}\ref{GenDiffExtendOrdinary} applied to the polynomial ring $W(B\{T\})$ over $A$, such a derivation indeed exists.\footnote{Notice that $W(B\{T\})$ naturally is a differential polynomial ring over $A$, but $\partial_{B\{T\}}^W$ is in general not the natural derivation of $W(B\{T\})$.} \par \medskip\noindent It remains to prove that there is a derivation $\partial_D^W$ of $W(D)$ as required. \par \smallskip\noindent \Claim 3 The ideal $I_D$ of $W(B\{T\})$ (see \ref{WeilD}) is a differential ideal for $\partial_{B\{T\}}^W$. \par \noindent \textit{Proof.} Let $f\in \ker(\pi_D)$. Then $W_{B\{T\}}(f)=\sum_{i=1}^\ell g_i\otimes b_i$, where $g_i=\lambda_i^{W(B\{T\})}(W_{B\{T\}}(f))$. By definition of $I_D$ it suffices to show that $\partial_{B\{T\}}^W(g_i)\in I_D$ \footnote{Notice that the module homomorphism $\lambda_i^{W(B(\{T\}))}:F(W(B\{T\}))\lra W(B\{T\})$ does not in general commute with the derivations.}. Now one checks that \[ W_{B\{T\}}(\partial_{B\{T\}}(f))=\sum_{i=1}^\ell\biggl(\partial_{B\{T\}}^W(g_i)+\sum_{j=1}^\ell \lambda_i(b_j)g_j\biggr)\otimes b_i \] Since $1\otimes b_1,\ldots,1\otimes b_\ell$ is a basis of $F(W(B\{T\}))$ over $W(B\{T\})$ we see that \[ \lambda_i^{W(B\{T\})}(W_{B\{T\}}(\partial_{B\{T\}}(f)))=\partial_{B\{T\}}^W(g_i)+\sum_{j=1}^\ell \lambda_i(b_j)g_j. \] The left hand side here is in $I_D$ by definition of $I_D$ and because $\ker(\pi)$ is differential for $\partial_{B\{T\}}$. As all $g_i\in I_D$ this entails $\partial_{B\{T\}}^W(g_i)\in I_D$. \hfill$\diamond$ \par \smallskip\noindent By Claim 3, the derivation $\partial_{B\{T\}}$ induces a derivation $\delta_D^W$ of $W(D)=W(B\{T\})/I_D$ such that $(W(D),\partial_D^W)$ is a differential $(A,d)$-algebra. It remains to show that $W_D$ is a differential $(B,\delta)$-algebra homomorphism, i.e., $W_D\circ \partial_D=(\partial_D^W\otimes \delta)\circ W_D$. This can be seen by a diagram chase as follows. Consider the diagram of maps \begin{center} \begin{tikzcd}[row sep=30pt, column sep=50pt] \ & D \arrow[dl, leftarrow, "\pi"'] \arrow[rr, "W_D"] \arrow[dd, "\partial_D" near end ] & \ & W(D)\otimes B \arrow[dl, leftarrow, "W(\pi)\otimes \id_B"'] \arrow[dd, "\partial^W_D\otimes \delta"]\\ {B\{T\}} \arrow[rr, crossing over, "W_{B\{T\}}" near end] \arrow[dd,"\partial_{B\{T\}}"'] & \ & W({B\{T\}})\otimes B \\ \ & D \arrow[dl, leftarrow, "\pi"'] \arrow[rr, "W_D"' near start] & \ & W(D)\otimes B \arrow[dl, leftarrow, "W(\pi)\otimes \id_B"] \\ {B\{T\}} \arrow[rr, "W_{B\{T\}}"'] &\ & W({B\{T\}})\otimes B \arrow[from=uu, crossing over, "\partial^W_{B\{T\}}\otimes \delta" near start]\\ \end{tikzcd} \end{center} The claim is that the back side of this cube is commutative. Now, all other sides of the cube are commutative squares, because \begin{itemize} \item Bottom and top of the cube are identical and commute as a property of the classical Weil descent. \item The front of the cube commutes as we know the \factname\ already for $(B\{T\},\partial_{B\{T\}})$. \item The square on the left hand side commutes by choice of $(B\{T\},\partial_{B\{T\}})$. \item The square on the right hand side commutes by applying base change to $B$ to the definition of $\partial_D^W$. \end{itemize} Since $\pi$ is surjective, we see that the back of the cube also commutes. This finishes the proof of existence and uniqueness of $\partial_W$. From Claim 2 we see that the definition of $\partial_{B\{T\}}^W$ only depends on $\partial_{B\{T\}}$ and not on $\delta$, because the structure map $B\lra D$ is differential. But then by construction of $\partial_{D}^W$ after Claim 3, $\partial_{D}^W$ only depends on $\partial_D$ and not on $\delta$. \end{proof} \begin{FACT}{Theorem}{liebracket} Let again $B$ be an $A$-algebra that is finitely generated and free as an $A$-module and let $D$ be a $B$-algebra. \par Let $\Der_B(D)$ be the set of all $\partial\in \Der(D)$ for which there are derivations $d$ of $A$ and $\delta$ of $B$ such that the structure maps of $B$ and $D$ are differential maps $(A,d)\lra (B,\delta)$ and $(B,\delta)\lra (D,\partial)$, respectively.\footnote{If the structure morphism of $D$ as a $B$-algebra is injective and we think of the structure maps $A\lra B$ and $B\lra D$ as inclusions, then $\Der_B(D)$ is the set of all derivations of $D$ that restrict to derivations on $A$ and $B$.} \par Then $\Der_B(D)$ is an $A$-submodule and a Lie subring of $\Der(D)$ and the map $\Der_B(D)\lra \Der(W(D))$ that sends $\partial$ to the derivation $\partial^W$ defined in \ref{partialDW}, is an $A$-module and a Lie ring homomorphism. Explicitly, given $\partial_1,\partial_2\in \Der_B(D)$ we have \begin{enumerate}[(i)] \item $(a_1\partial_1+a_2\partial_2)^W=a_1\partial_1^W+a_2\partial_2^W$ for all $a_1,a_2\in A$. \item $[\partial_1,\partial_2]^W=[\partial_1^W,\partial_2^W]$. In particular, $\partial_1^W,\partial_2^W$ commute if $\partial_1,\partial_2$ commute. \end{enumerate} \end{FACT} \begin{proof} In each case, the derivation of $W(D)$ on the right hand side turns $W(D)$ into a differential $A$-algebra, when $A$ is furnished with the derivation $a_1d_1+a_2d_2$ and $[d_1,d_2]$ respectively. By uniqueness in \ref{partialDW} we thus only need to verify the defining equation of the left hand side for the right hand side. \par \smallskip\noindent (i) Using \ref{GenDiff}\ref{GenDiffLieAndTensor}(a) we get $((a_1\partial_1^W+a_2\partial_2^W)\otimes (a_1\delta_1+a_2\delta_2))\circ W_D= a_1(\partial_1^W\otimes \delta_1)\circ W_D+a_2(\partial_2^W\otimes \delta_2)\circ W_D= a_1W_D\circ \partial_1+a_2W_D\circ \partial_2=W_D\circ (a_1\partial_1+a_2\partial_2)$, since $W_D$ is an $A$-algebra homomorphism. \par \smallskip\noindent (ii) Using \ref{GenDiff}\ref{GenDiffLieAndTensor}(b) we get $([\partial_1^W,\partial_2^W]\otimes [\delta_1,\delta_2])\circ W_D= [\partial_1^W\otimes \delta_1,\partial_2^W\otimes \delta_2]\circ W_D= (\partial_1^W\otimes \delta_1)\circ (\partial_2^W\otimes \delta_2)\circ W_D- (\partial_2^W\otimes \delta_2)\circ (\partial_1^W\otimes \delta_1)\circ W_D= (\partial_1^W\otimes \delta_1)\circ W_D\circ \partial_2- (\partial_2^W\otimes \delta_2)\circ W_D\circ \partial_1= W_D\circ \partial_1\circ \partial_2- W_D\circ \partial_2\circ \partial_1=W_D\circ [\partial_1,\partial_2]. $ \end{proof} \par \noindent Theorems \ref{partialDW} and \ref{liebracket} establish \begin{FACT}{:The differential Weil descent.}{DiffWeilDescent} Let $A$ be a ring and let $d=(d_i)_{i\in I}$ be a family of derivations of $A$. A differential $(A,d)$-algebra is an $A$-algebra $C$ together with derivations $(\eta_i)_{i\in I}$ of $C$ such that the structure map $A\lra C$ is a differential morphism $(A,d_i)\lra (C,\eta_i)$ for all $i\in I$. Let $(A,d)\text{-}\kat{Alg}$ be the category of differential $(A,d)$-algebras whose morphisms are ring homomorphisms preserving the appropriate derivations. \par We fix a differential $(A,d)$-algebra $(B,\delta)$, with $\delta=(\delta_i)_{i\in I}$, such that $B$ is finitely generated and free as an $A$-module. Then \begin{enumerate}[(i)] \item The functor $F^\mathrm{diff}:(A,d)\text{-}\kat{Alg}\lra (B,\delta)\text{-}\kat{Alg}$ that sends $(C,\eta)$ to $(C\otimes B,(\eta_i\otimes \delta_i)_{i\in I})$ has a left adjoint $W^\mathrm{diff}:(B,\delta)\text{-}\kat{Alg}\lra (A,d)\text{-}\kat{Alg}$, which we call the \notion{differential Weil descent} (or differential Weil restriction) from $(B,\delta)$ to $(A,d)$. It sends $(D,\partial)$ to $(W(D),\partial ^W)$ where $\partial^W=(\partial_i^W)_{i\in I}$ with $\partial_i^W$ as defined in \ref{partialDW}, and a morphism $f$ to $W(f)$. \item Let $(C,\eta)\in (A,d)\text{-}\kat{Alg}$ and let $(D,\partial)\in (B,\delta)\text{-}\kat{Alg}$. Then the bijection \[ \tau(D{,}C):\Hom_{A\text{-}\kat{Alg}}(W(D),C)\lra \Hom_{B\text{-}\kat{Alg}}(D,F(C)),\ \phi\longmapsto {F}(\phi)\circ {W}_D \] from the classical Weil descent \ref{WeilD} restricts to a bijection \[ \Hom_{(A,d)\text{-}\kat{Alg}}(W^\mathrm{diff}(D,\partial),(C,\eta))\lra \Hom_{(B,\delta)\text{-}\kat{Alg}}((D,\partial),F^\mathrm{diff}(C,\eta)). \] \item If $(D,\partial)\in (B,\delta)\text{-}\kat{Alg}$ and the derivations $\partial=(\partial_i)_{i\in I}$ are Lie commuting with structure coefficients $a_{ij}^k\in A$ $(i,j,k\in I)$, i.e., for fixed $i,j\in I$ only finitely many of the $a_{ij}^k$'s are nonzero and \[ [\partial_i,\partial_j]=\sum_{k\in I}a_{ij}^k\partial_k \quad (i,j\in I), \] then also the derivations $(\partial_i^W)_{i\in I}$ of $W(D)$ are Lie commuting with structure coefficients $a_{ij}^k$. \end{enumerate} \end{FACT} \begin{proof} By Theorem \ref{partialDW}, the map $W_D:D\to F(W(D))$ is differential. Hence, if a morphism $\phi:W(D)\to C$ is differential, so is ${F}(\phi)\circ {W}_D$. Thus the map $\tau(D,C)$ restricts to differential morphisms as claimed in (ii). Now recall from \ref{WeilD} (and \ref{LeftAdjoint}\ref{LeftAdjointDaggerImpliesLA}) that $W(f)$ is the unique map that corresponds to the morphism $W_{D'}\circ f:D\to F(W(D'))$ under the bijection $\tau(D,W(D'))$. As the latter morphism is differential, $W(f)$ must be differential. This entails (i), see \ref{LeftAdjoint}\ref{LeftAdjointLA}. Item (iii) follows immediately from \ref{liebracket}. \end{proof} \par \medskip \par Working over differential fields, Theorem \ref{DiffWeilDescent} has the following consequence at the level of rational points. This gives a geometric interpretation (in the sense of Kolchin's differential algebraic geometry \cite{Kolchi1973}) of the differential Weil descent. \begin{FACT}{Corollary}{consequence} Let $K$ be a field equipped with derivations $\delta=(\delta_i)_{i\in I}$. Suppose $L/K$ is a finite separable field extension. Recall that the derivations $(\delta_i)_{i\in I}$ extend uniquely from $K$ to $L$. Then, given a differential $L$-algebra $D$, by \ref{DiffWeilDescent}, there is a natural one-to-one correspondence between the differential $L$-points of $D$ and the differential $K$-points of $W^{\operatorname{diff}}(D)$. \end{FACT} \begin{FACT}{Remark}{remarkOnLeSMos} Assume $K$ and $L$ as in Corollary \ref{consequence}. Further assume that $\delta=\{\delta_1,\dots,\delta_m\}$ is a finite collection of commuting derivations and $D$ is the differential coordinate ring of an affine differential variety over $L$, say $D=L\{V\}$. In the case when $L$ has a $K$-basis $b_1,\dots,b_\ell$ of constants (meaning that $\delta_i(b_j)=0$ for all $i,j$), then a construction of the differential Weil descent $W^{\operatorname{diff}}(L\{V\})$ appears in \cite[\S5]{LeSMos2016}. However, a basis of constants does not always exist, as we point out in Example \ref{example1} below. \end{FACT} \begin{FACT}{Example}{example1} We work in the ordinary case $\Delta=\{\delta\}$. Let $K=\mathbb Q(t)$ with $\delta t=1$ and consider the finite extension $L=K(b)$ where $b^2=t$. Then the (unique) induced derivation on $L$ is given by $\delta b=\frac{1}{2b}=\frac{b}{2t}$. Fix the basis $\{1,b\}$ of $L$ as a $K$-module. Consider the differential variety $V$ given by $\delta x=0$ (i.e., $V$ is simply the constants of $\U$) viewed as a differential variety over $L$. The differential Weil descent $W^{\operatorname{diff}}(V)$ is obtained as follows; write $x$ as $x_1+x_2b$ and compute $$\delta(x_1+bx_2)=\delta x_1+(\delta b) x_2+b \delta x_2=\delta x_1+\frac{b}{2t} x_2+b \delta x_2=\delta x_1+\left(\frac{x_2}{2t}+\delta x_2\right)b.$$ Thus, $W^{\operatorname{diff}}(V)$ is the affine differential variety over $K$ given by the equations $$\delta x_1=0\quad \text{ and }\quad \delta x_2+\frac{x_2}{2t}=0.$$ Note that this is not contained in a product of the constants, as one might expect. Of course, if $\delta(b)$ were zero we would instead obtain the equations $\delta x_1=0$ and $\delta x_2=0$ (which would occur if $\delta$ were trivial on $K$, for instance). \end{FACT} \renewcommand\href[3][]{}  \par \end{document}